\documentclass[a4paper,12pt]{article}

\usepackage{amsmath}
\usepackage{amssymb}

\newtheorem{sth}{Theorem}
\newtheorem{sprop}[sth]{Proposition}
\newtheorem{slemma}[sth]{Lemma}

\newtheorem{scor}[sth]{Corollary}
\newcommand{\qed}{\,\,{\bf QED}}
\newcounter{std}

\begin{document}
\title{Gaps between classes of matrix monotone functions}
\author{Frank Hansen, Guoxing Ji and Jun Tomiyama \footnote{The second author was partially supported by the National Natural  Science Foundation of China(No.10071047) and the Excellent Young
  Teachers Program of MOE, P.R.C.}}

\date{April 16, 2002}

\maketitle
\section{Introduction} Almost seventy years have passed since K. L{\"o}wner \cite{L} proposed the notion of operator monotone functions. A real, continuous function 
$ f\colon I\to {\mathbf R} $ defined on an (non trivial) interval $I$ is said to be 
matrix monotone of order $ n $ if
\begin{equation}\label{definition of monotonicity}
x\le y\quad\Rightarrow\quad f(x)\le f(y)
\end{equation}
for any pair of self-adjoint $ n\times n $ matrices $x, y$ with eigenvalues in $I$. 
We denote by $ P_n(I) $ the set of such functions. A function $ f:I\to{\mathbf R} $ 
is said to be operator monotone if it is matrix monotone of arbitrary orders. We evidently
have $ P_{n+1}(I)\subseteq P_n(I) $ for each natural number $ n, $ and
$$
P(I)=\displaystyle\bigcap_{n=1}^\infty P_n(I)
$$
is the set of operator monotone functions defined on $ I. $ 
If (\ref{definition of monotonicity}) holds for any pair of self-adjoint elements 
$ x,y $ in a $ C^* $-algebra
$ A $ with spectra contained in  $ I, $ then we say that $ f $ is $ A $-monotone.

L{\"o}wner characterized the set of matrix monotone functions of order $ n $ in terms of positivity of certain
determinants (the so called L{\"o}wner determinants and the related Pick determinants) and proved that a function
is operator monotone if and only if it allows an analytic continuation to a Pick function,
that is an analytic function defined in the complex upper half plane with non-negative imaginary part.
Dobsch \cite{Do} continued L{\"o}wner's investigation and gave an alternative characterization of matrix monotonicity 
which we shall use in this paper.

Forty years after L{\"o}wner's work W. Donoghue published a comprehensive book on the subject in which
he refined Dobsch' necessary and sufficient condition for a function on an interval to be matrix monotone of order 
$n$ \cite[Chapter 7, Theorem VI and Chapter 8, Theorem V]{D}. Donoghue then
asserted \cite[p.~84]{D} that with this insight one may recognize that the classes $ P_n(I) $ are all distinct
for different values of $ n. $
We shall denote this as the (asserted) existence of gaps between the different classes of matrix monotone functions.

However, both L{\"o}wner's and Dobsch' conditions for matrix monotonicity of order $ n $ 
are very hard to check even 
for $ n=3, $ and explicit examples of functions showing such gaps are given by Donoghue 
only for $ n = 1 $ and $ n=2. $ Now another almost thirty years have passed after 
Donoghue's work and there are still, to our knowledge, no examples in the literature showing the gaps between 
$ P_n(I) $ and $ P_{n+1}(I) $ for arbitrary 
natural numbers $n.$ The purpose of this article is to prove exactly the existence 
of such gaps for every $ n. $ We  also characterize, for any natural number
$ n, $ the C*-algebras $ A $ with the property that any function 
$ f\in P_n(I) $ is $ A $-monotone.  
It is interesting to notice that this question is closely connected to the problem of matricial structure of operator algebras with respect to positive linear maps.

\section{The gap between $P_n(I)$ and $P_{n + 1}(I)$ }
For a positive integer $n$ let $g_n(t)$ be the polynomial defined by
\begin{equation}\label{g_n}
g_n(t) = t + \frac{1}{3}t^3 + \cdots + \frac{1}{2n - 1}t^{2n - 1}.
\end{equation}
Following the notations in \cite{D} we consider the matrix valued function associated with  $g_n(t)$ and given by
\[
M_n(g_n;t ) = \left(\frac{g^{(i + j - 1)}_n (t)}{(i + j - 1)!}\right)_{i,j=1}^n.
\]
The following lemma is an application of standard arguments from the theory
of moment problems for Hankel matrices.

\begin{slemma}\label{lemma} The matrix $M_n(g_n;0)$ is positive definite.
\end{slemma}
{\it Proof.}
We set
\[
b_k = \frac{1}{2}\int ^1_{-1} t^k dt\quad\mbox { for \(k = 0,1,2,\ldots\).}
\]
and calculate
\[
b_k = \left\{
\begin{array}{ll}
(k + 1)^{-1}\quad&\mbox{ if $k$ is even,}\\[1ex]
0 &\mbox{ if $k$ is odd.}
\end{array}\right.
\]
Hence, we can write $ g_n $ as
\[
g_n(t) = b_0 t + b_1 t^2 + \cdots + b_{2n - 2} t^{2n - 1}
\]
and therefore obtain
\[
g_n^{(i + j - 1)}(0) =  (i + j - 1)!\cdot b_{i + j - 2}\qquad i,j=1,\dots,n.
\]
Consequently
\[
M_n(g_n;0) = \bigl(b_{i + j - 2}\bigr)_{i,j=1}^n.
\]
Now take a vector \( c = (c_1,c_2,\ldots,c_n)\in{\mathbf C}^n\) and calculate
$$
\begin{array}{rl}
\bigl(M_n(g_n;0\bigr)c\mid c)& =
\displaystyle\sum^n_{i=1}\sum^n_{j=1} b_{i + j - 2}\, c_j\overline{c_i}\\[2ex]
&=\displaystyle\frac{1}{2} \sum^n_{i=1}\sum^n_{j=1}\int^1_{-1}
t^{i + j - 2}c_j\overline{c_i}dt\\[2ex]
&=\displaystyle\frac{1}{2}\int_{-1}^1\Bigl| \sum^n_{i=1}c_i t^{i-1}\Bigr|^2 dt.
\end{array}
$$
It follows that the matrix $M_n(g_n;0)$ is positive semidefinite. Moreover, if \(M_n(g_n;0)c = 0\) we see that
\[
\sum^n_{i=1}c_i t^{i-1} = 0\quad\quad\mbox{a.e.}
\]
Since this is a polynomial, it is identically zero on the interval $[- 1,1]$. All entries of the vector $c$ are therefore zero and the matrix $ M_n(g_n;0) $ is positive definite.
\qed\vskip 2ex

With this lemma we can show the existence of a gap between $P_n(I)$ and $P_{n+1}(I)$ for any positive integer $n$ and any nontrivial interval $I$ different from the whole real line.

\begin{sth} For any natural number $ n $ there exists a real number $ \alpha_n>0 $ and a function
$ g_n:[0,\alpha_n[\to{\mathbf R} $ such that
\begin{list}{(\arabic{std})}{\usecounter{std}\labelwidth=2em}

\item $ g_n $ is matrix monotone of order $ n $ on $ [0,\alpha_n[. $

\item $ g_n $ is not matrix monotone of order $ n+1 $ on $ [0,\alpha_n[, $ nor is it
matrix monotone of order $ n+1 $ on any subinterval.

\end{list}
\end{sth}
{\it Proof.} Consider the polynomial $g_n(t)$ introduced in the proof of Lemma~\ref{lemma}. By the continuous dependence of eigenvalues of matrices as a function of their entries, there exists by 
Lemma~\ref{lemma} a positive number $\alpha_n$ such that the matrix function $M_n(g_n;t)$ is positive definite 
for $ t\in [0,\alpha_n[. $ Since $ g^{(2n-3)} $ in addition is positive and convex on $ ]0,\alpha_n[ $ we conclude, 
cf. \cite[Chap.~8, Theorem V]{D}, that the
function $ g_n $ is matrix monotone of order $ n $ on the interval $ [0,\alpha_n[. $ 
The last principal matrix of order 3 of the matrix $M_{n+1}(g_n; t)$ is given by
\[
\left(
\begin{array}{ccc}
\displaystyle\frac{1}{2n - 3} + (n - 1)t^2&t &\displaystyle\frac{1}{2n - 1}\\
t & \displaystyle\frac{1}{2n - 1}& 0\\
\displaystyle\frac{1}{2n - 1} & 0 &0
\end{array}
\right )
\]
and this matrix has determinant $ -(2n-1)^{-3} $ regardless of the value of  $ t. $ The matrix
$M_{n+1}(g_n;t)$ is thus not positive semi-definite and the function $g_n$ is not matrix monotone of order $n + 1$ on 
any subinterval $ J\subseteq [0,\alpha_n[. $ This completes the proof.
\qed\vskip 2ex

Consider the concrete function $ g_n $ defined in equation~\ref{g_n}. A calculation shows that the largest
possible value of $ \alpha_2 $ is $ 1. $
It is exceedingly difficult to calculate the largest possible value for $ n\ge 3. $ 

\begin{sprop}
Let either $ I=[a,b[ $ or $ I=[a,\infty[ $ for real numbers $ a<b $ and take $ \alpha>0. $
Then there exists a bijection $ h:[0,\alpha[\to I $ such that both $ h $ and the inverse map 
are operator monotone. Likewise, with $ J=]a,b] $ or $ J=]-\infty,b], $ there exists a 
bijection $ g:]0,\alpha]\to J $ such that both $ g $ and the inverse map are operator monotone. 
\end{sprop}
{\it Proof.} An affine map of the form $ t\to ct+d $ with $ c>0 $ is operator monotone, and so is the inverse
map $ t\to c^{-1}(t-d). $ We may therefore assume that $ \alpha=1, $ $ I=[0,1[ $ or $ I=[0,\infty[, $ and
$ J=]0,1] $ or $ J=]-\infty,0]. $ The function
$$
h(t)=t(1+t)^{-1}\quad\mbox{with inverse}\quad h^{-1}(t)=t(1-t)^{-1}
$$
is a bijection of $ [0,\infty[ $ to $ [0,1[. $ Likewise is the function
$$
g(t)=(1-t)^{-1}\quad\mbox{with inverse}\quad h^{-1}(t)=1-t^{-1}
$$
a bijection of $ ]\infty,0] $ to $ ]0,1]. $ The assertion follows since $ h, $ $ h^{-1}, $ $ g, $ $ g^{-1} $ are
all operator monotone functions, cf.~\cite{BS, H2}.
\qed\vskip 2ex

Notice that we cannot find a bijection $ h:]0,1[\to{\mathbf R} $ such that both $ h $ and $ h^{-1} $ are
operator monotone. An operator monotone function defined on the whole real line is necessarily affine, cf.~\cite{D}. Its
range is therefore either a constant or the whole real line.

\begin{scor}\label{cor}
Let $ I=[a,b[ $ or $ I=[a,\infty[ $ for real numbers $ a<b. $ For any natural number $ n $ there exists a function
$ f_n:I\to{\mathbf R} $ such that

\begin{list}{(\arabic{std})}{\usecounter{std}\labelwidth=2em}

\item $ f_n $ is matrix monotone of order $ n $ on $ I. $

\item $ f_n $ is not matrix monotone of order $ n+1 $ on $ I, $ nor is it
matrix monotone of order $ n+1 $ on any subinterval.

\end{list}

\end{scor}

Let $ I $ be any open real interval and take $ t_0\in I. $ Bendat and Sherman proved in \cite[Theorem 3.2]{BS} that a function $ f:I\to{\mathbf R} $ is matrix convex of order $ n, $ if and only if the function
$$
F(t)=\frac{f(t)-f(t_0)}{t-t_0}
$$
is matrix monotone of order $ n. $ Notice that $f,$ for $ n\ge 2, $ automatically is differentiable and 
$ F(t_0)=f'(t_0). $ One may set $ F(t_0)=(f(t_0)_+ + f(t_0)_-)/2 $ for $ n=1. $ We also have $ f(t)=f(t_0)+F(t)(t-t_0). $ 

\begin{scor}
Let $ I=[a,b[ $ or $ I=[a,\infty[ $ for real numbers $ a<b. $ For any natural number $ n $ there exists a function
$ f_n:I\to{\mathbf R} $ such that

\begin{list}{(\arabic{std})}{\usecounter{std}\labelwidth=2em}

\item $ f_n $ is matrix convex of order $ n $ on $ I. $

\item $ f_n $ is not matrix convex of order $ n+1 $ on $ I, $ nor is it
matrix convex of order $ n+1 $ on any subinterval.

\end{list}

\end{scor}

The statement follows by combining Bendat and Sherman's result with Corollary \ref{cor}.

\section{Characterization of $C^*$-algebras in terms of matrix monotone functions}

As we have discussed in \cite{GT}, we may regard the question of monotonicity of functions as a kind of nonlinear version of 
the problem of matricial structure of operator algebras. Recall that a positive linear map $\tau$ from a $C^*$-algebra $A$ to a $C^*$-algebra $B$ is said to be $n$-positive if the map 
\[
\tau_n : \bigl(a_{ij}\bigr)_{i,j=1}^n \to \bigl(\tau(a_{ij})\bigr)_{i,j=1}^n
\]
is a positive map from $M_n(A)$ to $M_n(B)$. If $\tau$ is $n$-positive for all positive integers, then
it is said to be completely positive.

Although the introduction of these notions by Stinespring \cite{S} is of a much later date than the work of L{\"o}wner,
they have turned out to be very important notions for the matricial structure of operator algebras i.e.
$C^*$-algebras and von Neumann algebras. One may simply recognize this aspect by the recent publication \cite{ER} by Effros and Ruan. Meanwhile examples of $n$-positive maps which are not $(n + 1)$-positive had been investigated, and it had been discussed for
which types of $C^*$-algebras $A$ every $n$-positive map from or to $A$ for an another $C^*$-algebra $B$ is also 
$(n + 1)$-positive. In this sense, gaps between $P_{n+1}(I)$ and $ P_n(I) $ are nonlinear versions of the above sort of problems. 
We are thus naturally led to the problem of the characterization of those $C^*$-algebras $ A $ on which every matrix monotone functions of order $n$ is $ A $-monotone. The following theorem is a generalization of a previous result \cite[Theorem 1]{GT} where the two last authors essentially treated the gap between $P_1(I)$ and $P_2(I)$. 
In this investigation we reach the same kind of $C^*$-algebras as in the study of positive linear maps 
by the third author \cite {T1}.

\begin{sth} Let $ A $ be a C*-algebra, and let $ I $ be an interval of the form $ I=[a,b[ $ or $ I=[0,\infty[ $ for real numbers $ a<b. $ The following assertions are equivalent:

\begin{list}{(\arabic{std})}{\usecounter{std}\labelwidth=2em}

\item Every matrix monotone function of order $ n $ defined on $I$ is $A$-monotone.

\item The dimension of every irreducible representation of $A$ is less or equal to $n.$

\item Every $n$-positive linear map from/to $A$ for another C*-algebra $B$ is completely positive.

\end{list}
\end{sth}
{\it Proof.} $(1)\Rightarrow (2)$:\, We first notice that we, without loss of generality, may choose the interval
$ I=[0,\infty[. $
Suppose that $A$ had an irreducible representation $\pi$ on a Hilbert space $H$ whose dimension is greater than $n$. Take an $(n + 1)$-dimensional projecton $e$ in $H$. We then have 
\(\pi(A)e = B(H)e\) by \cite[Theorem 4.18]{Ta}, hence
\[
e\pi(A)e = eB(H)e = B(eH) \cong M_{n+1}.
\]
Let $B$ be the $C^*$-subalgebra of $A$ defined by setting
$$
B=\{a\in A\mid \pi(a)e=e\pi(a)\}.
$$
By the theorem cited above, the restriction of $\pi$ to $B$ is a $*$-homomorphism onto $eB(H)e$. 
We choose a function $f$ in $P_n(I)$ which is not matrix monotone of order 
$n + 1,$ cf. Corollary~\ref{cor}.  Let $c$ and $d$ be arbitrary positive elements in $eB(H)e$ with $ c \leq d. $
It is easily verified that we can find positive elements $a$ and $b$ in $B$ such that 
$a \leq b$, \(\pi(a) = c\) and \(\pi(b) = d.\) Since $a\leq b$ we obtain $f(a) \leq f(b) $ by the assumption, whence
\[
f(c) = f(\pi(a))= \pi (f(a)) \leq \pi(f(b)) = f(\pi(b)) = f(d).
\]
Therefore, $f$ is matrix monotone of order $n + 1$ on $I,$ a contradiction.

\noindent $(2)\Rightarrow (1):$\, Take a function $f$ in $P_n(I)$ and let $a$ and $b$ be self-adjoint elements in $A$ with spectra contained in $I$ such that $a \leq b.$ We consider an irreducible representation $ \pi $ of $A.$ Since
also the spectra of $\pi(a)$ and $\pi(b)$ are contained in $I,$ we obtain by the assumptions that
\[
\pi(f(a)) = f(\pi(a)) \leq f(\pi(b)) = \pi(f(a)).
\]
It follows that $f(a) \leq f(b),$ thus $f$ is $A$-monotone.

\noindent $(2)\Leftrightarrow (3): $\, The assertion is proved in \cite{T1}.\qed

\nocite{H}

{\small
\noindent Frank Hansen (frank.hansen@econ.ku.dk)\\
Institute of Economics, Copenhagen University, Denmark.\\[1ex]
Guoxing Ji (gxji@snnu.edu.cn)\\
College of Mathematics and Information Science, Shaanxi Normal University, Xi'an, 710062, P.R. China.\\[1ex]
Jun Tomiyama (jtomiyama@fc.jwu.ac.jp)\\
Department of Mathematics and Physics, Japan Women's University, Mejirodai Bunkyo-ku, Tokyo, Japan.}

\end{document}